# An On-Line Algorithm for Improving Performance in Navigation

**Avrim Blum***
School of Computer Science
Carnegie Mellon University
Pittsburgh, PA 15213
avrim@theory.cs.cmu.edu

**Prasad Chalasani**
School of Computer Science
Carnegie Mellon University
Pittsburgh, PA 15213
chal@cs.cmu.edu

## Abstract

Recent papers have shown optimally-competitive on-line strategies for a robot traveling from a point $s$ to a point $t$ in certain unknown geometric environments. We consider the question: Having gained some partial information about the scene on its first trip from $s$ to $t$, can the robot improve its performance on subsequent trips it might make? This is a type of on-line problem where a strategy must exploit *partial information* about the future (e.g., about obstacles that lie ahead). For scenes with axis-parallel rectangular obstacles where the Euclidean distance between $s$ and $t$ is $n$, we present a deterministic algorithm whose *average* trip length after $k$ trips, $k \leq n$, is $O(\sqrt{n/k})$ times the length of the shortest $s$-$t$ path in the scene. We also show that this is the best a deterministic strategy can do. This algorithm can be thought of as performing an optimal tradeoff between search effort and the goodness of the path found. We improve this algorithm so that for *every* $i \leq n$, the robot's $i$th trip length is $O(\sqrt{n/i})$ times the shortest $s$-$t$ path length. A key idea of the paper is that a *tree* structure can be defined in the scene, where the nodes are portions of certain obstacles and the edges are "short" paths from a node to its children. The core of our algorithms is an on-line strategy for traversing this tree optimally.

## 1 Introduction

Imagine you have just moved to a new city; you are at your home and must travel to your office, but you do not have a map. Let's assume you know your coordinates and those of your office. A collection of papers in recent literature have studied on-line competitive strategies for quickly traveling from point A to point B for problems of this sort. But now, suppose you have reached your office, spent the day there, and it is time to go home. You could retrace your path, but you now have some information about the city (what you saw on your way to work in the morning) and would like to do better. The next morning you have even more information and so on. What is a good strategy that allows your performance at each stage to be as good as possible, and to improve with experience? Perhaps you might even design your paths explicitly so as to gain more information for future trips. This is the sort of problem we consider here.

Specifically, we consider the scenario (examined in [10, 3, 7]) where the start point $s$ and target $t$ are in a 2-dimensional plane filled with non-overlapping axis-parallel rectangular obstacles. A point robot begins at $s$, and knows its current position and that of the target, but it does *not* know the positions and extents of the obstacles; it only finds out of their existence as it encounters them. In the problem considered in previous papers, the robot must travel from $s$ to $t$ circumventing the obstacles. We call this the *one-trip* problem. In this paper we consider a robot that may be asked to make *multiple trips*, going back and forth between $s$ and $t$.

It is intuitive that in the multiple trips problem, information from previous trips must be used if one hopes to improve the path in later trips. In particular, on later trips some but not all of the obstacles in the scene are known. We therefore have a type of on-line problem that is different from the standard scenario, in that *partial information* about the future (e.g., about obstacles that lie ahead) must be exploited to achieve good performance.

A particular arrangement of $s$, $t$, and the obstacles is called the *scene*. Let $n$ denote the *Euclidean* distance between $s$ and $t$ in the scene, where the obstacles have width and height at least 1. Papadimitriou and Yannakakis [10] showed for the one-trip problem a lower bound of $\Omega(\sqrt{n})$ on the ratio of the distance traveled by the robot to the actual shortest path (the *competitive ratio*) for any deterministic algorithm. Blum, Raghavan, and Schieber [3] describe an algorithm whose performance matches this bound. Whether randomization can help improve upon this bound is an open question, although a lower bound was obtained by [7].

For the multiple trips problem, there are several ways one might formalize the intuitive goal described in the first paragraph. One natural way is to consider the total distance traveled in the $k$ trips between $s$ and $t$ and examine the ratio of this to $k$ times the shortest path. Thus, for $k = 1$, the previous results give a ratio of $\Theta(\sqrt{n})$. For $k = n$ it is not hard to see that one can achieve a ratio of $O(1)$ by simply performing a search of cost proportional to $n$ times the shortest path length, to find the shortest path on the first trip. Our main result is to show an optimal smooth transition. For $k \leq n$ we present an algorithm whose competitive ratio is $O(\sqrt{n/k})$ and give an $\Omega(\sqrt{n/k})$ lower bound for deterministic algorithms. A key idea of the algorithm is to optimally traverse a certain *tree* structure based on the obstacles in the scene.

*Supported in part by an NSF Postdoctoral Fellowship.

Notice that the "cumulative" formulation allows one to search "hard" on the first trip to find a short path, and then use this short path on the remaining $k-1$ trips, which is what we do. This result can be thought of as one that shows how to optimally trade off exploration effort with the goodness of the path found. We in addition show how to modify the algorithm so that on the $i$th trip, its ratio *for that trip* is $O(\sqrt{n/i})$. So, this algorithm is optimal for each prefix cumulative cost and in addition does not spend too much effort on any one trip. Thus, the algorithm can be viewed as one that optimally improves its performance with each trip, achieving the intuitive goal described at the start of this paper.

**Related Work.** In the machine learning literature (especially reinforcement learning), some authors have addressed problems similar to the multiple-trips problem [4, 12, 8]. The problem of efficiently visiting several destinations has been examined in the robotics literature [9] but not from the viewpoint of competitive analysis. A variety of models and algorithms for efficient, complete exploration of an unknown environment (rather than just visiting particular destinations) have been studied by previous authors [2, 5, 6, 11].

## 2 The model, and some preliminaries

Let $\mathcal{S}(n)$ denote the class of scenes where the Euclidean distance between $s$ and $t$ is $n$. We define $s$ to be at the origin $(0,0)$. As mentioned above, we assume that the width and height of the obstacles are at least 1 (this in essence defines the units of $n$) and for simplicity that the $x$-coordinates of the corners of obstacles are integral. Thus no more than $n$ obstacles can be placed side by side between $s$ and $t$. We assume that when obstacles touch, the point robot can "squeeze" between them.

To simplify the exposition, we take $t$ to be the infinite vertical line (a "wall") $x = n$ and require the robot only to get to any point on this line; this is the Wall Problem of [3]. Our algorithms can be easily extended to the case where $t$ is a point, using the "Room Problem" algorithms of [3] or [1]. This modification is sketched in the appendix.

We assume that the only sensing available to the robot is *tactile*, that is, it discovers an obstacle only when it "bumps" into it. It will be convenient to assume that when the robot hits a rectangle, it is told which corner of the rectangle is nearest to it, and how far that corner is from its current position. As in [3], our algorithms can be modified to work without this assumption with only a constant factor penalty. We do not describe these modifications in this paper.

Consider a robot strategy $R$ for making $k$ trips between $s$ and $t$. Let $R_i(S)$ be the distance traveled by the robot in the $i$th trip, in scene $S$. Let $L(S)$ be the length of the shortest obstacle-free path in the scene between $s$ and $t$. We define the *cumulative* $k$-trip competitive ratio as $\rho(R,n,k) = \max_{S\in\mathcal{S}(n)} \frac{R^{(k)}(S)}{kL(S)}$, where $R^{(k)}(S) = \sum_{i=1}^{k} R_i(S)$ is the *total* distance traveled by the robot in $k$ trips. That is, $\rho(R,n,k)$ is the ratio between the robot's total distance traveled in $k$ trips, and the best $k$-trip distance. We define the *per-trip* competitive ratio for the $i$th trip as $\rho_i(R,n) = \max_{S\in\mathcal{S}(n)} \frac{R_i(S)}{L(S)}$.

Our main results are the following. First, we show for any $k$, $n$, and deterministic algorithm $R$, that $\rho(R,n,k) = \Omega(\sqrt{n/k})$. Second, we describe a deterministic algorithm that given $k \leq n$ achieves $\rho(R,n,k) = O(\sqrt{n/k})$. Finally, we show an improvement to that algorithm that in any scene achieves $\rho_i(R,n) = O(\sqrt{n/i})$ for *all* $i \leq n$.

## 3 A Lower Bound for $k$ Trips

**Theorem 1 ($k$-trip Cumulative Lower Bound)** *For $k \leq n$, the ratio $\rho(R,n,k)$ is at least $\Omega(\sqrt{n/k})$, for any deterministic algorithm $R$.*

**Proof:** Since $R$ is deterministic, an adversary can simulate it and place obstacles in $\mathcal{S}$ as follows. Recall that $s$ is the point $(0,0)$.

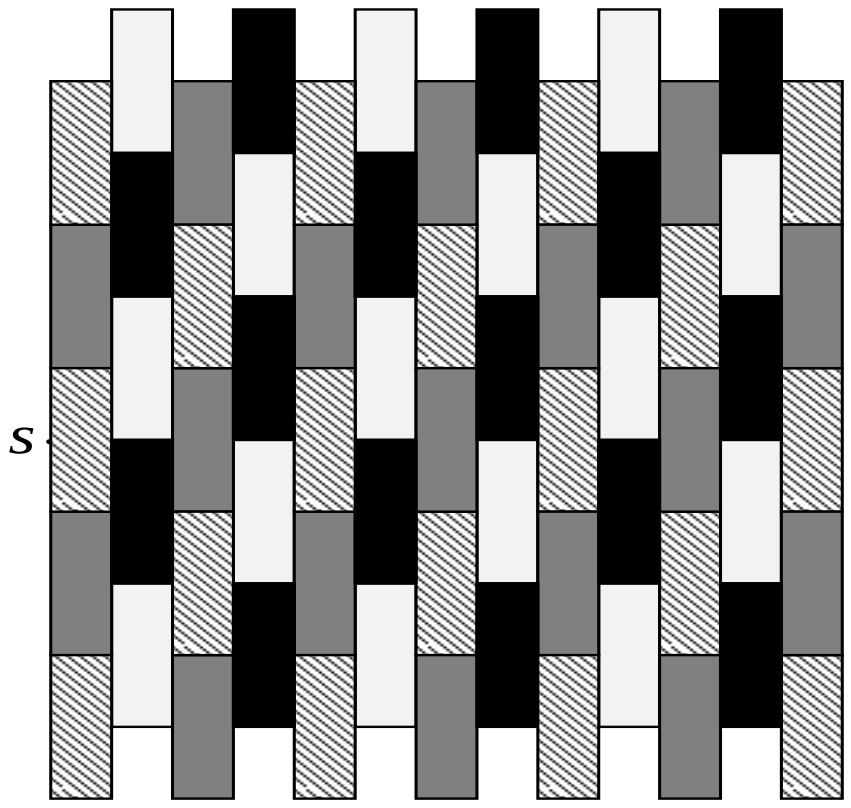

The adversary first places obstacles of fixed height $h \geq \sqrt{n}$ and width 1, in a full "brick pattern" on the entire plane, as shown above, with $s$ at the center of the left-side of an obstacle. (Recall that the point robot can "squeeze" between bricks). The adversary simulates $R$ on this scene, notes which obstacles it has touched at the end of $k$ trips, then removes all other obstacles from the scene. This is the final scene that the adversary creates for the algorithm, and say it contains $M$ obstacles. The brick pattern ensures that $R$ must have hit at least one brick at every integer $x$-coordinate, so $M \geq n$. Further, this arrangement forces the robot to hit a brick at every integer $x$-coordinate on *every* trip. Whenever it hits a brick, it must move vertically up or down a distance $h/2$, so its total $k$-trip distance $R^{(k)}$ is at least $nkh/2$.

We now show that there is a path from $s$ to the wall of length at most $O(\sqrt{R^{(k)}h})$. Imagine the full brick pattern to be built out of four kinds of bricks (red, blue, yellow and green, say) arranged in a periodic pattern as shown in the above figure. This arrangement has the following property: for each color, to go from a point on an obstacle of that color to a point on any other of the same color, the robot must move a distance at least $h/2$. Out of the $M$ obstacles hit by the robot, at least $M/4$ must have the same color, say blue. So regardless of how the robot moved, since it has visited $M/4$ blue obstacles, we have $R^{(k)} \geq Mh/8$, which implies $M \leq 8R^{(k)}/h$.

We claim there is a non-negative integer $j \leq \sqrt{M}$ such that at most $\sqrt{M}$ obstacles have centers at the $y$-coordinate $jh$. This is because a given obstacle intersects at most one $y$-coordinate of the form $jh$, and there are $M$ obstacles. Thus, there is a path to $t$ that goes vertically to the $y$-coordinate $jh$, then horizontally along this $y$-coordinate, going around at most $\sqrt{M}$ obstacles. The total length of this path is at most $h\sqrt{M}+h\sqrt{M}+n$, which is at most $3h\sqrt{M}$ since $n \leq M$ and $\sqrt{n} \leq h$. Since $M \leq 8R^{(k)}/h$, this path is in fact of length at most $3\sqrt{8hR^{(k)}}$. Thus the $k$-trip ratio is at least $R^{(k)}/(3k\sqrt{8hR^{(k)}})$. Recalling that $R^{(k)} \geq nkh/2$, this is at least $\frac{1}{12}\sqrt{n/k} = \Omega(\sqrt{n/k})$. ∎

It is not hard to see that this lower bound also holds for the case where $t$ is a point rather than a wall.

# 4 An Optimal Algorithm

The second and main result in this paper is a deterministic algorithm for making $k$ trips that achieves a cumulative ratio of $O(\sqrt{n/k})$, matching the lower bound proved above. To keep the discussion simple, we assume for now that the algorithm (robot) knows both the length $L$ of the shortest path from $s$ to $t$, and the number of trips $k$. We later show how these assumptions can be removed.

Our approach is to devote the first trip to searching "hard" for a short path, and then to use this short path on the remaining $k-1$ trips. In particular, the algorithm first performs an "exploratory" walk of length $O(L\sqrt{nk})$, which has the property that an $s$-$t$ path of length $O(L\sqrt{n/k})$ can be composed out of portions of the walk. In other words, the robot travels a distance of only $O(\sqrt{k})$ times the $L\sqrt{n}$ bound for the one-trip problem, and finds an $s$-$t$ path that is guaranteed to be $\Omega(\sqrt{k})$ times shorter than $L\sqrt{n}$. The algorithm then uses this path on the remaining $k-1$ trips, achieving a cumulative $k$-trip ratio of $O(\frac{1}{kL}(L\sqrt{nk} + kL\sqrt{n/k})) = O(\sqrt{n/k})$.

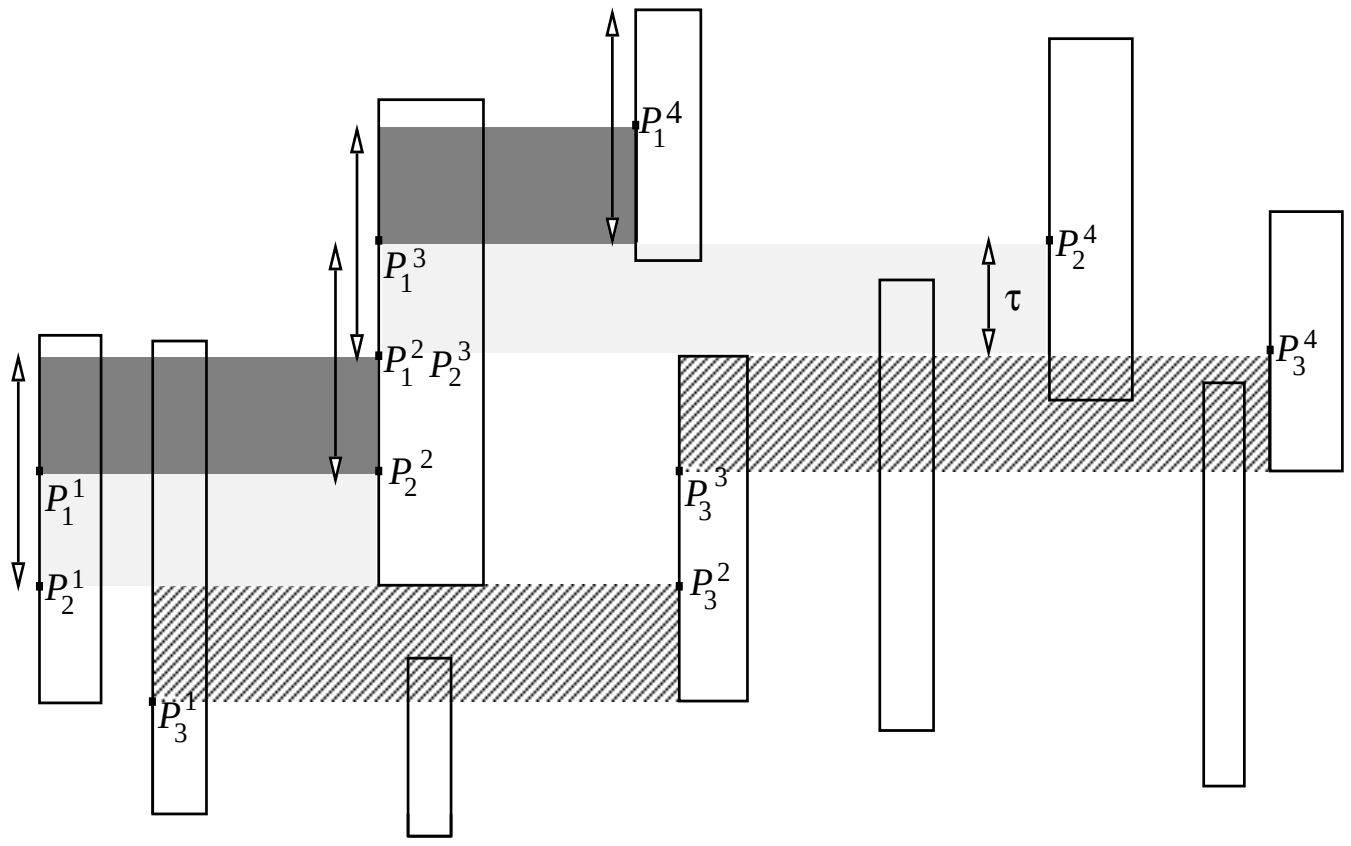

Figure 1: A collection of 3 disjoint fences with 4 posts each. The solid rectangles are the obstacles. The bands of different fences are shaded differently.

## 4.1 Fences and Fence-posts

We first introduce some terms that will be useful to picture the working of the algorithm. We will use the words up, down, left, and right to mean the directions $+y, -y, -x, +x$ respectively. When we say point $A$ is above, below, behind, or ahead of a point $B$ we will mean that $A$ is in the $+y, -y, -x, +x$ direction respectively from $B$. Finally, vertical (horizontal) motion is parallel to the $y$ (respectively, $x$) axis.

In a given scene, we define a *$\tau$-fence* $F$ in terms of *fence-posts* as follows.[1] A *$\tau$-post* is a (vertical) portion of height $2\tau$ of the left-edge of an obstacle. When we say a post is at a point $P$, we mean its center is at $P$, and we will often identify a post with its center point. A *$\tau$-fence* $F$ is a sequence of $\tau$-posts at points $P^1 = (X^1, Y^1), P^2 = (X^2, Y^2), \ldots, P^M = (X^M, Y^M)$ such that for $m = 1, 2, \ldots, M-1$:

$$X^m \leq X^{m+1} \quad (1)$$
$$Y^{m+1} = Y^m + \tau \ \left(= Y^1 + (m-1)\tau\right) \quad (2)$$

We use a subscript to distinguish between different fences. For example, the $m$'th post of fence $F_i$ is denoted by $P_i^m$ and its coordinates are $(X_i^m, Y_i^m)$. In Fig. 1, the sequence of $\tau$-posts $\langle P_1^1, P_1^2, P_1^3, P_1^4 \rangle$ form a $\tau$-fence $F_1$.

The axis-parallel rectangular region of height $\tau$ whose opposite corners are the centers of $P^{m-1}$ and $P^m$ is called a *band* $B^m$. Note that the inequality (1) is not strict, so two consecutive posts $P^{m-1}$, $P^m$ may lie along the same obstacle, so that the band $B^m$ is empty (i.e., has zero area). Two fences are said to be *disjoint* if their *non-empty* bands do not overlap. Thus the three fences in Fig. 1 are disjoint.

Given a fence with $M$ posts, if we consider all points in the plane between the lines $y = Y^1$ and $y = Y^M$

[1] What we call a *fence* is similar to what was called a *sweep* in [3].

*excluding* those in the bands, then these points form two regions, one on each side of the fence. Any path that goes from one region to the other without going above $y = Y^M$ or below $y = Y^1$ is said to *cross* the fence. *The most important property of a $\tau$-fence is that it costs at least $\tau$ to cross*, since to cross the fence, one must cross one of its bands, and thus travel a vertical distance of at least $\tau$. Since the non-empty bands of disjoint fences do not overlap, it follows that a collection of $k$ disjoint $\tau$-fences costs at least $k\tau$ to cross.

### 4.2 A High-Level View of the Search Algorithm

Using the notion of fences, we can give a high-level view of the algorithm for walking a distance $O(L\sqrt{nk})$ and finding a path of length $O(L\sqrt{n/k})$.

Recall, the start point $s$ is $(0,0)$ and the wall $t$ is at $x = n$. The algorithm maintains a rectangular *window* of height $2L$ bounded by the lines $x = 0$, $x = n$, $y = L$, $y = -L$. Notice that the shortest $s$-$t$ path must lie inside the window. The value of the threshold $\tau$ is chosen to be $L/\sqrt{nk}$. The robot repeatedly executes the following two steps. If at any time during these steps it hits the wall, then the algorithm halts.

1. Go along a "greedy" right-down path (i.e., move right until an obstacle is hit, then down to the corner of the obstacle, and repeat) to the bottom of the window. This path has length $O(L+\Delta x)$ where $\Delta x$ is $x$-distance between the initial and final positions in this step. Call this path a *group-transition* path.
2. Walk in such a way as to discover a group of $k$ disjoint $\tau$-fences $F_1, F_2, \ldots, F_k$, with each fence extending from the bottom of the window to the top of the window. Suppose the net $x$-motion is $\Delta x$. We must perform the walk such that: (a) the total distance traveled is $O(kL + k\tau\ \Delta x)$, *and* (b) there is a path discovered of length $O(L + \tau\ \Delta x)$ that crosses the entire collection of fences, i.e., connects the first post of $F_1$ to the last point on the walk, which is the last post of fence $F_k$. Performing a walk that achieves both (a) and (b) is the main difficult step of the algorithm. Call the path of (b) a *group-crossing path.*

Because a fence costs $\tau = L/\sqrt{nk}$ to cross, there can be at most $\sqrt{nk}$ disjoint $\tau$-fences in the window and therefore the algorithm will repeat the above steps at most $\sqrt{n/k}$ times. Note that if $L$ is not known, we can just guess a value, and if the above steps repeat more than $\sqrt{n/k}$ times we can double our guess and repeat the entire procedure. Thus there is only a constant factor penalty for not knowing $L$.

Since the $x$-motions do not overlap between groups of fences found in step (2), the $\Delta x$ terms add to at most $n$, so the total distance traveled is $O(kL\sqrt{n/k} + nk\tau) = O(L\sqrt{nk})$. (The greedy paths are a low order term). In addition, a group can be crossed by a group-crossing path of length $O(L+\tau\ \Delta x)$, and a group-transition path of length $O(L + \Delta x)$ connects its end to the beginning of the next group-crossing path. Thus, again since the $\Delta x$ terms sum to at most $n$, and there are at most $\sqrt{n/k}$ groups, there is a path from $s$ to $t$ composed of alternating group-transition and group-crossing paths of total length $O(L\sqrt{n/k} + n\tau) = O(L\sqrt{n/k})$.

The remainder of the paper describes how step 2 above is done. Unlike the strategy in [3], we will not actually traverse each fence in order. Instead, the algorithm works by finding a collection of fences whose posts can be thought of as nodes on a *tree.* The root of this tree is the first post of the first fence ($P_1^1$), and an edge is a "short" path from a node to its child. As a byproduct of traversing the edges of this tree we will not only have "cheaply" found the desired collection of fences, but the tree path from the root to the last post of the last fence will be a "cheap" group-crossing path. The next section describes this tree structure.

### 4.3 The Fence-Tree

We first introduce some useful terminology. A $\tau - path$ is the path of the robot when it moves to the right along some line $y = y_0$ as follows: If it hits an obstacle whose nearest corner is less than $\tau$ away, it goes around that corner to the point on the opposite side with $y$-coordinate $y_0$ and continues to the right; if it hits either a $\tau$-post or the wall, it stops. E.g., in Fig. 2, the path from $A$ to $\tau$-post $P_1^2$ is a $\tau$-path. Observe that a $\tau$-path has vertical motion at most $2\tau$ at every (integer) $x$-coordinate on the path, so:

**Fact 1** *A $\tau$-path between two points $(x, y)$ and $(x + \delta x, y)$ has length at most $\delta x + 2\tau\ \delta x$.*

Given a $\tau$-post $P_0$ in a scene $\mathcal{S}$, the unique $k \times M$ *$\tau$-fence-tree with root $P_0$* is a binary tree defined as follows, when it exists. This tree has $kM$ nodes, where each node is a certain $\tau$-post in the scene and is denoted by $P_i^m$ for $i = 1, 2, \ldots, k$ and $m = 1, 2, \ldots, M$. (The reason for this notation is that the posts form a disjoint collection of $k$ fences with $M$ posts each, but let us ignore the fence interpretation for now.) In the following definition, when we say $\tau$-post $P'$ is an *up-child* (*down-child*) of $\tau$-post $P$ we mean that $P'$ is the $\tau$-post at the end of the $\tau$-path starting from the *top* (*bottom*) of post $P$. If this $\tau$-path hits the wall, then this child does not exist and the desired tree does not exist in the scene. So if $P'$ is a child of $P$ then there is a path from $P$ to $P'$ consisting of a vertical portion of length $\tau$, and a $\tau$-path. We call this path an *up-edge* or a *down-edge* depending on whether $P'$ is an up-child or down-child of $P$. We are now ready to define the tree.

*The post $P_1^1$ is the root of this tree and is at the given initial post $P_0$.* The locations of the remaining posts (nodes) $P_i^m$ of the tree are specified by the following

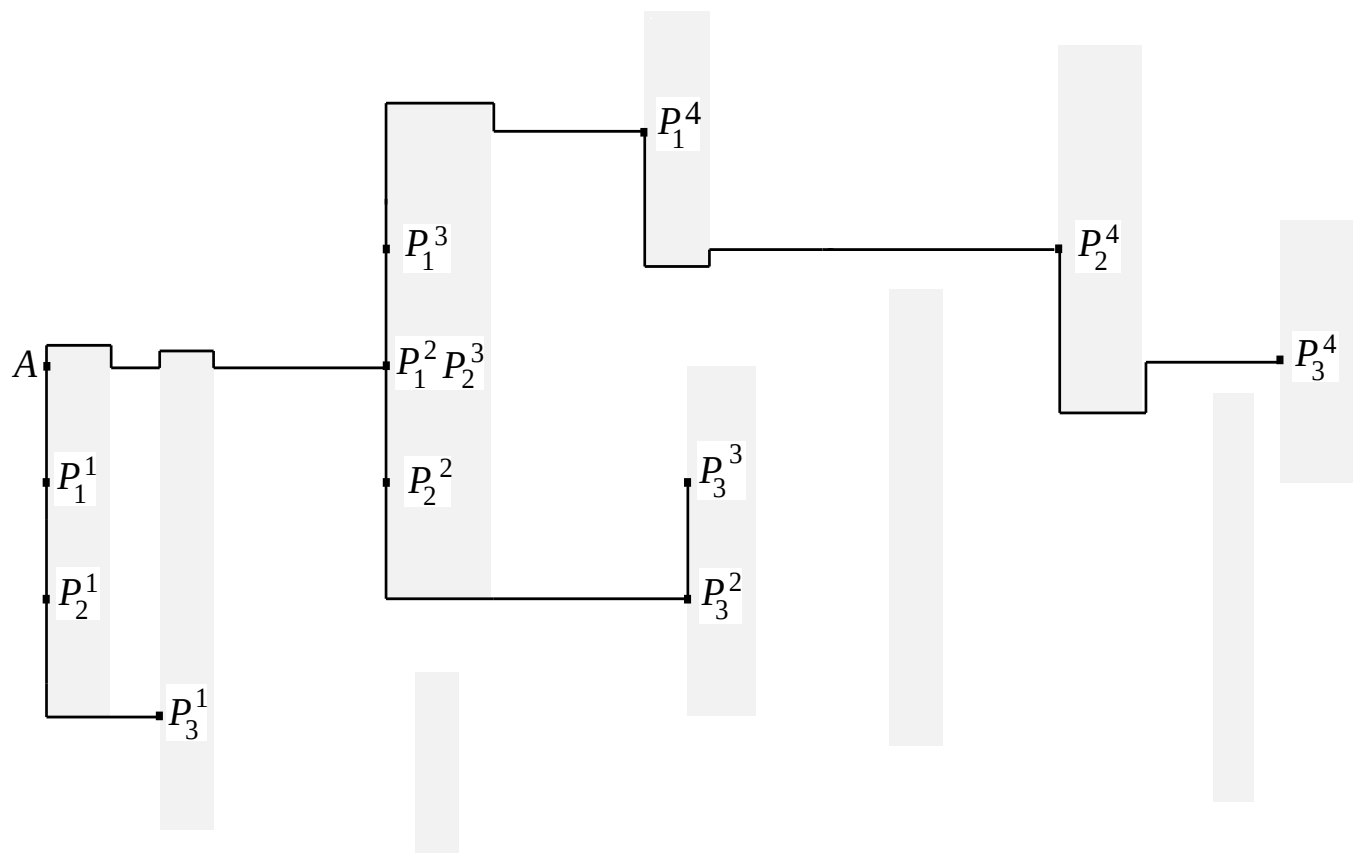

Figure 2: A $3 \times 4$ $\tau$-fence-tree. The shaded rectangles are the obstacles, and solid lines are tree edges. The fences corresponding to this tree are shown in Fig. 1

rules. Examples of these rules can be seen in Fig. 2. Recall that the coordinates of $P_i^m$ are $(X_i^m, Y_i^m)$.

1. For each $m \in \{2, \ldots, M\}$, $P_1^m$ is the up-child of $P_1^{m-1}$. Thus $Y_1^m = Y_1^1 + (m-1)\tau$.
2. For each $i \in \{2, \ldots, k\}$, $P_i^1$ is the down-child of $P_{i-1}^1$. Thus $Y_i^1 = Y_1^1 - (i-1)\tau$.
3. For each $i \in \{2, \ldots, k\}$ and $m \in \{2, \ldots, M\}$, $P_i^m$ is either the down-child of $P_{i-1}^m$ or the up-child of $P_i^{m-1}$. Which will be the case is decided as follows:
   (a) If $X_{i-1}^m > X_i^{m-1}$ then $P_i^m$ is the down-child of $P_{i-1}^m$ (and $P_i^{m-1}$ does not have an up-child in the tree). E.g., $P_3^4$ is the down-child of $P_2^4$ but not the up-child of $P_3^3$ in Fig. 2.
   (b) Otherwise (i.e., $X_{i-1}^m \leq X_i^{m-1}$) $P_i^m$ is the up-child of $P_i^{m-1}$ (and $P_{i-1}^m$ does not have a down-child in the tree.) E.g., $P_3^3$ is the up-child of $P_3^2$ but not the down-child of $P_2^3$ in Fig. 2.

It is easy to see that rule (3) above implies the following two facts:

$$X_i^m \geq \max\{X_{i-1}^m,\ X_i^{m-1}\}, \tag{3}$$

$$Y_i^m = Y_i^1 + (m-1)\tau = Y_1^1 + (m-i)\tau. \tag{4}$$

One implication of these facts is that for each $i$, the posts $P_i^1, \ldots, P_i^M$ form a $\tau$-fence $F_i$, and the $m$'th post of a fence $F_i$ is *exactly* $\tau$ *higher* than the $m$'th post of $F_{i+1}$. We thus say that the fence $F_i$ is *above* $F_{i+1}$. Also, from relation (3), the $m$'th post of $F_i$ is never to the right of the $m$'th post of $F_{i+1}$, and it is not hard to see that this ensures that the fences defined by this tree are disjoint. *Thus a $k \times M$ $\tau$-fence-tree defines a collection of $k$ disjoint $\tau$-fences with $M$ posts each.* Fig. 1 shows the fences defined by the fence-tree of Fig. 2.

Notice that if the $y$-distance between $P_1^1$ and $P_k^M$ in the $k \times M$ $\tau$-fence-tree is $2L$, (and by equation (4) this implies $M = 2L/\tau + k$.) then the collection of $k$ fences satisfies our requirements for step 2 in the high-level idea of the last section, namely: (a) All the fences extend across a window of height $2L$, and (b) The path in the tree from $P_1^1$ to $P_k^M$ is the "group-crossing" path that crosses all the $k$ fences. The theorem below states that this group-crossing path is "short": i.e., has length $O(L + \tau\Delta x)$, and that the total length of all the edges is of the same order as our desired bound on the cost of finding the $k$ fences, namely $O(kL + k\tau\Delta x)$.

**Theorem 2** *Suppose there is a $k \times M$ $\tau$-fence-tree consisting of fences $F_1, F_2, \ldots, F_k$, with $Y_k^M - Y_1^1 = 2L$ and $X_k^M - X_1^1 = \Delta x$, such that $L \geq n$, $\tau = L/\sqrt{nk}$, and $k \leq n$. Then:*
*(a) The unique path in the tree from $P_1^1$ to $P_k^M$ has length at most $4L + 3\tau\Delta x$;*
*(b) The total length of all the edges in the fence-tree is at most $k(3L + 3\tau\Delta x)$*
**Proof.** (See Appendix)

It remains to show how to traverse this fence-tree efficiently, i.e., with cost $O(kL + k\tau\Delta x)$. Note that we cannot traverse the tree depth-first since in general a node $P_i^m$ can be located only after *both* its possible parents $P_i^{m-1}$ and $P_{i-1}^m$ have been identified (rule 3 in the fence-tree definition). Thus, it may happen that we discover $P_i^{m-1}$ first, then $P_{i-1}^m$ and we find by rule 3 that $P_i^m$ is the up-child of $P_i^{m-1}$, so we must re-visit $P_i^{m-1}$ in order to find $P_i^m$. A naive traversal strategy might suffer a high cost when doing this type of re-visiting. For example, we might try to find the posts of the tree fence-by-fence, starting with fence $F_1$. A worst-case tree for this strategy is one where each post $P_i^m$ of $F_i, i > 1$, is the down-child of $P_{i-1}^m$. Thus if we have just discovered the post $P_i^m$, $i > 1$, then in order to find $P_i^{m+1}$ we must re-visit $P_{i-1}^{m+1}$, which may require us to follow tree edges all the way up to $F_1$, and back down to $F_{i-1}$. In fact, scenes can be constructed where this strategy will have a cost of $\Omega(k^2 L)$.

The next section shows our algorithm for finding the fence tree of Theorem 2 with a walk of length $O(kL + k\tau\Delta x)$. From the lower bound on the cumulative $k$-trip ratio, it follows that this algorithm is optimal up to a constant factor.

## 4.4 An Optimal Algorithm to Traverse a Fence Tree

The fence-tree traversal problem is this: $L, \tau, k$ are given, such that $L \geq n$, $k \leq n$ and $\tau = L/\sqrt{nk}$. Let $M = 2L/\tau + k = 2\sqrt{nk} + k$. Initially the robot is at a $\tau$-post at $(x_0, -L)$. If a $k \times M$ $\tau$-fence tree exists with root $P_1^1 = (x_0, -L)$, then the robot must traverse all edges in this tree, and eventually arrive at the post $P_k^M$. If such a tree doesn't exist, then the robot must arrive at the wall. In either case, if the $x$-distance between the initial and final positions is $\Delta x$, the total distance moved by the robot should not exceed a constant times $(kL + k\tau\Delta x)$. We claim that the procedure FindFenceTree (described below) accomplishes just this.

The procedure maintains two variables to keep track of the "progress" made so far on each fence $F_i$: $M_i$ denotes the number of posts found so far in $F_i$ (the "$y$" progress); $(X_i, Y_i)$ denotes the coordinates of the last (rightmost, i.e., most recently found) post on $F_i$ ($X_i$ is the "$x$ progress"). Recall that we assume the robot is initially at a $\tau$-post at $(x_0, -L)$. So initially, we set $(X_1, Y_1) = (x_0, -L)$, $M_1 = 1$, and all the other $M_i$ and $X_i$ are set to 0. The initial $Y_i$ values are not relevant. The procedure repeatedly tries to *extend* some fence, i.e., find the next post of the fence. The variable $i$ keeps track of Recall that if post $P_i^m, i > 1, m > 1$, has been found then both $P_{i-1}^m$ and $P_i^{m-1}$ must already have been found, so at any stage, a fence has at least as many posts as the one below:

$$M_1 \geq M_2 \geq \ldots M_k. \tag{5}$$

The procedure maintains the following two loop invariants, which trivially hold initially.

**$x$-ordering.** All fences below the current fence $F_i$ are $x$-ordered, in that $X_{i+1} \leq X_{i+2} \leq \cdots \leq X_k$,

**Almost-$x$-ordering.** At most one post (i.e., the last one) of a fence is strictly to the right of the last post of any lower fence. Thus, for each $j$: $X_j^{M_j - 1} \leq X_\ell$ for all $\ell > j$.

Invariant [Almost-$x$-ordering] implies the following fact. If at some stage the next post of a fence $F_i$ is the down-child of a post $P$ of the fence $F_{i-1}$ above, then $P$ is ahead of the last post of $F_i$ (rule 3a in the fence tree definition) and therefore $P$ is the *last* post of $F_{i-1}$. Thus a new edge in the tree can only come from the last post of a fence – a property that does not hold for the naive algorithm we described before.

**Fence-extension rules.** If $P_i^m$ is the last post of a fence $F_i$, $i > 1$ (so $M_i = m$) then the following three conditions determine how $P_i^{m+1}$ can be found:

1. Either (a) $M_{i-1} = m + 1$ and $X_{i-1} \leq X_i$, or (b) $M_{i-1} > m+1$. In case (a), clearly $X_{i-1}^{m+1} = X_{i-1} \leq X_i = X_i^m$, and in case (b), by invariant [Almost-$x$-ordering], $X_{i-1}^{m+1} \leq X_i^m$. Thus by rule 3(b) in the fence-tree definition, $P_i^{m+1}$ is the up-child of $P_i^m$. In this case we say that $F_i$ *can be extended from* $F_i$.
2. For some fence $F_j$ above $F_i$, $X_j > X_i$ and $M_j = m+1$. From relation 3 we know that $X_j = X_j^{m+1} \leq X_{j+1}^{m+1} \leq \ldots \leq X_{i-1}^{m+1}$, so $X_{i-1}^{m+1} > X_i = X_i^m$ and therefore by rule 3(a) in the fence-tree definition, $P_i^{m+1}$ is the down-child of $P_{i-1}^{m+1}$ (which may not have been found). In this case we say that *fence* $F_i$ *can be extended from* $F_{i-1}$. Note that in this case we also know that when $P_i^{m+1}$ is found eventually, $X_i^{m+1} \geq X_j$ must hold.
3. If neither of the above conditions hold, then we cannot as yet determine how $F_i$ can be extended. Note that in this case, the following must hold: $X_{i-1} \leq X_i$ and $M_{i-1} = M_i$.

We now describe the procedure FindFenceTree in detail. The procedure uses two subroutines GoDown $(i, M_i, M_{i+1})$ and GoBackDown $(i, M_i, M_{i+1})$ to move from the last post of fence $F_i$ to the last post of fence $F_{i+1}$. The first subroutine is used when $X_i \leq X_{i+1}$, and the second is used when $X_i > X_{i+1}$. When using these routines the robot may move along paths that are not edges of the fence-tree. This is necessary since it is possible to construct a scene where *any* traversal strategy (even a randomized one) that only moves on the edges must pay at least $\Omega(k^2 L)$. We prove the correctness and bound the costs of these subroutines in the appendix, so here we assume that they have the desired effect. The current coordinates of the robot at any time are denoted $(x, y)$. We say a fence $F_i$ is *ahead* of fence $F_j$ if $X_i > X_j$; $F_i$ is *behind* $F_j$ if $X_i < X_j$. Recall that variable $i$ is the "current fence" number, which is initially 1. The robot is initially at the starting $\tau$-post $P_1^1 = (x_0, -L)$.

**Procedure** FindFenceTree.

The procedure repeatedly checks the conditions of the following 6 cases in sequence, and executes the action corresponding to the *first* case that applies. It will be clear that the invariant [$x$-ordering] is maintained since (a) whenever the current fence $F_i$ is ahead of the fence $F_{i+1}$ below, the robot goes down to work on $F_{i+1}$ (cases 1, 2) and (b) the robot returns to the fence above (case 6) only when the current fence $F_i$ is not ahead of $F_{i+1}$. Also, invariant [Almost-$x$-ordering] is maintained because of invariant [$x$-ordering] and the fact that a post is added (in cases 2, 3) to a fence $F_j$ only when $F_j$ is *not* ahead of $F_{j+1}$. Note that the last case always applies if none of the earlier ones apply.

1. *The current fence $F_i$ is ahead of the next lower fence $F_{i+1}$, and $F_{i+1}$ can be extended from $F_{i+1}$* (i.e., $i < k$, $X_i > X_{i+1}$, and $M_i > M_{i+1} + 1$).

   Go down to the last post of $F_{i+1}$ using procedure GoBackDown $(i, M_i, M_{i+1})$; $i \leftarrow i + 1$.

2. *The current fence $F_i$ is ahead of the next lower fence $F_{i+1}$* (i.e., $i < k$ and $X_i > X_{i+1}$.)

   Traverse a new down-edge, i.e., go vertically down a distance $\tau$, then on a right $\tau$-path until a $\tau$-post. If the wall is reached then HALT.

$M_{i+1} \leftarrow M_{i+1} + 1$; $(X_{i+1}, Y_{i+1}) \leftarrow (x, y)$;
$i \leftarrow i + 1$.
(If case 1 did not apply, it must be that $M_i = M_{i+1} + 1$, so $F_{i+1}$ can be extended from $F_i$.)

3. *The current fence $F_i$ can be extended from $F_i$.* (i.e., either:
(a) $i = 1$ and $M_i < M$, or
(b) $i > 1$ and either
(i) $M_{i-1} > M_i + 1$, or
(ii)$M_{i-1} = M_i + 1$ and $X_{i-1} \leq X_i$.)

   Traverse a new up-edge, i.e., go vertically up a distance $\tau$, then right on a $\tau$-path to a $\tau$-post. If the wall is reached, then HALT;
$M_i \leftarrow M_i + 1$; $(X_i, Y_i) \leftarrow (x, y)$.

4. *The next lower fence $F_{i+1}$ can be extended from $F_{i+1}$, and either (a) the current fence is completed, or (b) it is determined by fence-extension rule (2) that the next post of $F_i$ will be ahead of the last post of $F_{i+1}$.* (i.e., $i < k$, $M_i > M_{i+1}$, and either (a) $M_i = M$ or (b) for some fence $F_u$ above $F_i$, $M_u = M_i + 1$ and $X_u > X_{i+1}$)

   Go down to the last post of $F_{i+1}$ using procedure GoDown $(i, M_i, M_{i+1})$; $i \leftarrow i + 1$.

5. *The last post of the lowest fence has been found.* (i.e., $i = k$ and $M_i = M$).

   In this case, HALT.

6. Return to the last post of $F_{i-1}$, i.e., to $(X_{i-1}, Y_{i-1})$, by the last path that was used to get to the current point; $i \leftarrow i - 1$.
(In this case it must be either that $F_i$ can be extended from the higher fence $F_{i-1}$, or that it is not yet possible to determine how $F_i$ can be extended (extension rule 3))

We show the correctness of this procedure:

**Theorem 3** *The procedure* FindFenceTree *terminates, and if the robot is not at the wall on termination, it will be at post $P_k^M$ and will have found all nodes and edges of the unique $k \times M$ $\tau$-fence-tree with the given root post $P_1^1$.*

**Proof:** On each iteration, the value of the fence-number $i$ either remains the same (case 3), increases by 1 (cases 1, 2, 4), or decreases by 1 (case 6). It is easy to see that $i$ is always at least 1 and at most $k$, and none of the values $M_j$ for $j = 1, 2, \ldots, k$ exceeds $M$. Whenever the value of $i$ remains the same, a new post is added to the fence $F_i$. In addition, whenever the value of $i$ increases by 1, a new post is added to $F_{i+1}$ in either the current or the subsequent iteration. Thus whenever $i$ remains the same or increases by 1, some $M_j$ increases by 1. Since $i$ is bounded above and below, and the $M_j$ are bounded above, the procedure must terminate. If the robot is not at the wall on termination, then case 5 must have applied, and the robot be at post $P_k^M$ and therefore have found all posts and edges on the tree. ■

In the appendix we show that the total distance walked by the robot in procedure FindFenceTree is at most $(60kL + 62k\tau\Delta x)$. This gives our main result:

**Theorem 4** *There is a deterministic algorithm $R$ for a robot that for any $k \leq n$ achieves $\rho(R, n, k) = O(\sqrt{n/k})$.*

## 4.5 An incremental algorithm

We describe here an improvement of our cumulative algorithm, so that the per-trip ratio on the $i$'th trip, for all $i \leq n$, is $O(\sqrt{n/i})$. Let us for simplicity say that we know $L$. From the earlier results in this paper, we know that by searching a distance at most $cL\sqrt{nk}$ we can find an $s$-$t$ path of length at most $c'L\sqrt{n/k}$, for some constants $c, c'$ and any $k \leq n$. Let us suppose that at the end of $i$ trips we know an $s$-$t$-path $\pi$ of length at most $c'L\sqrt{n/i}$. Then we know how to search with cost at most $cL\sqrt{n2i}$ and find a path of length at most $c'L\sqrt{n/2i}$. Let us denote by $\Pi$ the path we would have traveled if we did this entire search in one trip. In order to maintain a per-trip ratio of $O(\sqrt{n/i})$, we spread the work of $\Pi$ over the next $i$ trips as follows. Each trip consists of two phases: The first is a *search* phase, where we walk an additional portion of $\Pi$ of length $\frac{1}{i}cL\sqrt{n2i} = cL\sqrt{2n/i}$, starting from where we left off on the previous trip. We can always do this because the fences are in a tree structure, so that the last point in $\Pi$ during the previous search can always be reached from the start point by a known short path whose length adds only a small constant factor to the total trip length. Once the search phase is completed, we "give up" and enter the *follow* phase, where we complete the trip by joining (by a greedy path) the known path $\pi$ of length $c'L\sqrt{n/i}$, and following it to $t$. Thus our trip length is still $O(L\sqrt{n/i})$. Since in each such search-follow trip we traverse a portion of $\Pi$ of length $cL\sqrt{2n/i}$, and the length of $\Pi$ is at most $cL\sqrt{2ni}$, after $i$ trips we will have completely walked the path $\Pi$. So after the first $2i$ trips we have a path of length at most $c'L\sqrt{n/2i}$. We can then repeat this procedure. Since initially, we can find a path of length $c'L\sqrt{n}$ (to start off the induction), we have the following theorem:

**Theorem 5** *There is a deterministic algorithm $R$ that achieves for every $i \leq n$, a per-trip ratio on the $i$'th trip, $\rho_i(R, n)$, of $O(\sqrt{n/i})$.*

# References

[1] E. Bar-Eli, P. Berman, A. Fiat, and P. Yan. On-line navigation in a room. In *Proc. 3rd ACM-SIAM SODA*, 1992.

[2] M. Betke, R. Rivest, and M. Singh. Piecemeal learning of an unknown environment. In *Proc. 6th ACM Conf. on Computational Learning Theory*, pages 277–286, 1993.

[3] A. Blum, P. Raghavan, and B. Schieber. Navigating in unfamiliar geometric terrain. In *Proc. 23rd ACM STOC*, 1991.

[4] P. C. Chen. Improving path planning with learning. In *Proc. 9th Int'l Workshop on Machine Learning*, 1992.

[5] X. Deng, T. Kameda, and C. Papadimitriou. How to learn an unknown environment. In *Proc. 32nd IEEE FOCS*, pages 298–303, 1991.

[6] X. Deng and C. Papadimitriou. Exploring an unknown graph. In *Proc. 31st IEEE FOCS*, pages 355–361, 1990.

[7] H. Karloff, Y. Rabani, and Y. Ravid. Lower bounds for randomized $k$-server and motion-planning algorithms. In *Proc. 23rd ACM STOC*, pages 278–288, 1991.

[8] S. Koenig and R. G. Simmons. Complexity analysis of real-time reinforcement learning. In *Proc. AAAI*, pages 99–105, 1993.

[9] B. Oomen, S. Iyengar, N. Rao, and R. Kashyap. Robot navigation in unknown terrains using learned visibility graphs. *IEEE J. Robotics and Automation*, 3:672–681, 1987.

[10] C. Papadimitriou and M. Yannakakis. Shortest paths without a map. In *Proc. 16th ICALP*, 1989.

[11] R. L. Rivest and R. E. Schapire. Inference of finite automata using homing sequences. In *Proc. 21st ACM STOC*, pages 411–420, 1989.

[12] S. B. Thrun. Efficient exploration in reinforcement learning. Technical Report CMU-CS-92-102, Carnegie Mellon University, 1992.

# APPENDIX

For the proof of Theorem 2, Lemmas 7 and 8, we need the following facts. From equation (4), $2L = Y_k^M - Y_1^1 = (M-k)\tau$, so $M = 2L/\tau + k$. Since $L \geq n$, we have $\tau = L/\sqrt{nk} \geq 1$. From $k \leq n$ it follows that $k\tau = kL/\sqrt{nk} \leq L$. This implies $M\tau = 2L + k\tau \leq 3L$.

## Proof of Theorem 2

Part (a). There are exactly $(M+k-2)$ edges in the tree path from $P_1^1$ to $P_k^M$. Since the vertical portion of each edge has length $\tau$, and the $\tau$-path portions of the edges do not overlap in the $x$-direction, the total length of these edges is at most $(M+k-2)\tau + 2\tau\Delta x + \Delta x$, which is at most $(4L + 3\tau\Delta x)$ from the inequalities above.

Part (b). Note that we can associate each edge with a unique post, namely the one at the right-end of the edge. For any given post $P_i^m$ other than $P_1^1$, the $x$-distance to its parent is *at most* the $x$-distance $\delta x$ to its predecessor $P_i^{m-1}$ on the same fence. So the edge associated with this post has length at most $(\tau + 2\tau\delta x + \delta x)$. The sum of the $\delta x$ terms over all posts of the fence $F_i$ is the $x$-distance between the first and last posts of $F_i$, which is at most $\Delta x$. So the total length of the edges associated with the $M$ posts of a fence is at most $(M\tau + 2\tau\Delta x + \Delta x)$, which sums to $k(M\tau + 2\tau\Delta x + \Delta x)$ for $k$ fences. This last expression is at most $k(3L + 3\tau\Delta x)$ from the previous inequalities. ■

## Modification for Point-to-Point Navigation.

Our algorithms can be extended to the case where $t$ is a point rather than a wall, with essentially the same bounds, as follows. Let us assume for simplicity that the shortest path length $L$ is known. As before, if we do not know $L$, we can use the standard "guessing and doubling" approach and suffer only a constant factor penalty in performance. On the first trip, the robot can get to $t$ using the optimal point-to-point algorithms of [3] or [1], with a single-trip ratio of $O(\sqrt{n})$. Once at $t$, the robot creates a greedy up-left path and a greedy down-left path from $t$, within a window of height $4L$ centered at $t$. Note that the highest post in a $k \times M$ $\tau$-fence-tree is $M\tau \leq 3L$ above the root (which is always distance $L$ below $t$) and the lowest post is $k\tau \leq L$ below the root. So the robot is guaranteed to stay within a window of height $4L$ centered at $t$. Thus after the first trip, these greedy paths play the role of a wall; once the robot hits one of these paths, it can reach $t$ with an additional cost that is only a low-order term in the total cost.

## Procedures GoDown and GoBackDown and Cost Analysis

We first introduce some useful terminology. A *monotone right-up* path is a path that only goes right or up. Other monotone paths are defined similarly. A greedy path is thus a special kind of monotone path. Clearly, the length of a monotone path between $(x,y)$ and $(x',y')$ is $|x-x'| + |y-y'|$.

We now describe the procedures GoDown and GoBackDown used by the procedure FindFenceTree.

**Procedure** GoDown $(i, m, q)$
**Initial Conditions:**
*The robot is at a point* $(x,y)$ *satisfying*
$X_i^m \leq x \leq X_{i+1}^q$ *and* $Y_i^m \geq y \geq Y_{i+1}^q$.
$m \geq q$; $X_i^m \leq X_{i+1}^q$;
**Final Condition:**
*The robot is at* $P_{i+1}^q$.
**Begin**

1. *Go greedy down-left until* $y \leq Y_{i+1}^q + \tau$;

2. *Go greedy right-down until either*

   (a) *the post $P_{i+1}^q$ is reached, in which case:* **return** *, or*

   (b) *an edge is hit;*

3. *Follow the edges to the right until $P_{i+1}^q$ is reached.*

**End**

**Procedure** GoBackDown $(i, m, q)$
**Initial Conditions:**
*The robot is at $P_i^m$;*
$m \geq q+2,\ X_i^m > X_{i+1}^q \geq X_i^{m-1}$;
**Final Condition:**
*The robot is at $P_{i+1}^q$.*
**Begin**

1. *Follow tree path backward until* $x = X_{i+1}^q$;

2. *Go greedy down-left until either*

   (a) *the current point is in the rectangular region whose opposite corners are the centers of posts $P_j^{m-1}$ and $P_{j+1}^{m-1}$, for some $j < i$. In this case do the following until $j = i$:*
   GoDown $(j, m-1, m-1)$; $j \leftarrow j+1$,
   *or:*

   (b) $y = Y_i^{m-1}$,

3. GoDown $(i, m-1, q)$;

**End**

In the proofs ahead, we will need the following important property of the procedure FindFenceTree.

**Lemma 6** *For $1 < i \leq k$, when the procedure* FindFenceTree *moves up the robot from fence $F_i$ to fence $F_{i-1}$ (case 6), then*
*(a) If $X_i \geq X_{i-1}$ then $M_i = M_{i-1}$, and*
*(b) For every $s \geq i$, if $X_s < X_{i-1}$ then $M_s = M_{i-1} - 1$.*

**Proof:** Let us refer to the current iteration as "iteration B", and let $M_{i-1} = m$ in this iteration. Since case 3 does not apply in this iteration, if $X_i \geq X_{i-1}$ then $M_i = m$, and $M_i = m-1$ otherwise. This proves that (a) holds, and also that (b) holds for $s = i$.

To prove (b) we show that if for some $s \geq i$, $X_{s+1} < X_{i-1}$ and $M_s = m-1$, then $M_{s+1} = m-1$. Consider the *last* iteration (call it iteration A) before iteration B when the robot was at fence $F_s$. Thus case 4 cannot apply in iteration $A$. Also, no new posts can be found between iterations A and B. So if $X_{s+1} < X_{i-1}$ then $F_{s+1}$ must have as many posts as $F_s$, which is $m-1$. ■

We now prove the correctness and bound the cost of procedure GoDown.

**Lemma 7** *(1) Suppose $m \geq q$, $X_i^m \leq X_{i+1}^q$, and the robot is at a point $(x, y)$ such that $X_i^m \leq x \leq X_{i+1}^q$ and $Y_i^m \geq y \geq Y_{i+1}^q$. If the procedure* GoDown $(i, m, q)$ *is now invoked then the final position of the robot is at post $P_{i+1}^q$;*
*(2) The total distance walked during all calls to* GoDown *in case 4 of procedure* FindFenceTree *is at most $k(9L + 9\tau\Delta x)$, where $\Delta x = X_k^M - X_1^1$.*

**Proof:** Note that whenever the procedure Let us consider the path of the robot in each of the steps of the procedure GoDown $(i, m, q)$.

In step 1, the robot moves greedy down-left to a point $A$ with $y$-coordinate $Y_{i+1}^q + \tau$, which is also the $y$-coordinate of the bottom of $P_i^{q+1}$. (Note that if $m = q$ there is no motion in this step since $y \leq Y_{i+1}^q + \tau = Y_i^m$ is already true.) This path is bounded on the left by the posts of $F_i$, so the worst case step 1 path length is ($\tau+$ the length of a monotone path from $P_i^{q+1}$ to $P_i^m$).

In step 2, the robot moves greedy right-down from $A$ until it hits either (a) the post $P_{i+1}^q$ or (b) an edge of the tree. To see that one of these must occur, note two facts. First, recall that $A$ has the same $y$-coordinate as $P_{i+1}^q$. Second, $A$ cannot be to the left of $P_i^{q+1}$, and $F_{i+1}$ has only $q$ posts, so any posts of $F_i$ on the unique tree-path from the root to $P_{i+1}^q$ must be below $A$. Thus in the worst case the step 2 path goes vertically down from $A$ to some point $B$ on this tree path leading to $P_{i+1}^q$. The step 3 path is simply the tree path from $B$ to $P_{i+1}^q$. Since $P_{i+1}^q$ is exactly $\tau$ lower than $A$, the step 2 path is at most $\tau$ longer than the step 3 path. It is clear at this point that (1) holds, i.e., the final position of the robot is at post $P_{i+1}^q$.

To bound the total cost of this procedure, we show that in two *successive* calls GoDown $(i, m, q)$ and GoDown $(i, m', q')$ from the *same* fence $F_i$ in case 4 of procedure FindFenceTree, (a) the left end of the step 1 path in the second call is at least as high as and strictly to the right of the right end of the corresponding path in the first call, and (b) the tree-paths followed in step 3 in the two calls are distinct. Fact (a) implies that the total step 1 cost in calls to this procedure from fence $F_i$ is at most ($\tau\Delta x+$ the length of a monotone path from $P_i^1$ to $P_i^M$ ), which is $(\tau\Delta x + \Delta x + M\tau) \leq 3L + 2\tau\Delta x$. This sums to at most $k(3L + 2\tau\Delta x)$ for all $k$ fences. Combined with the fact that in different calls to this procedure from fence $F_i$, the step 3 path follows only edges associated with $F_{i+1}$, fact (b) implies that the total cost of this step in all calls is at most the total length of all tree edges, which is at most $k(3L + 3\tau\Delta x)$ (Theorem 2). This in turn implies that the total step 2 cost is at most $k\tau\Delta x + k(3L + 3\tau\Delta x)$. Thus the total cost of this procedure is at most $k(9L + 9\tau\Delta x)$ which proves (2).

To argue facts (a) and (b), we make two observations. If the call GoDown $(i, m, q)$ is made in case 4, $M_i = m$, $m > q$, and it is known that the *next* post of $F_i$ will be to the right of $P_{i+1}^q$,i.e., $X_i^{m+1} > X_{i+1}^q$ – this is the first observation. Also, when the robot returns to fence $F_i$ we know by lemma 6 that $M_{i+1}$ must be at least $m$, so in the *next* call (if any) GoDown $(i, m', q')$ we must have $q' \geq m$ – this is the second observation.

The left end of the step 1 path in the *next* (if any) call GoDown $(i, m, q)$ in the worst case is the bottom of $P_i^{q'+1}$, which is $P_i^{m+1}$ or some higher post on $F_i$. Since $X_i^{m+1} > X_{i+1}^q$ and we already know that $X_i^m \leq X_{i+1}^q$, it follows that $P_i^{m+1}$ is strictly to the right of $P_i^m$, which is the right end of the step 1 path in the first call. This implies fact (a). The right end of the step 3 path in the first call is at post $P_{i+1}^q$. The left end of the corresponding path in the next call (if any) has $x$-coordinate at least $X_i^{q'+1} \geq X_i^{m+1} > X_{i+1}^q$, from the two observations above. This implies fact (b). ■

We show the correctness and bound the cost of procedure GoBackDown.

**Lemma 8** *(1) If procedure* GoBackDown $(i, m, q)$ *is invoked (by procedure* FindFenceTree*) then on termination the robot will be at* $P_{i+1}^q$.
*(2) The total distance walked by the robot in* **all** *calls to this procedure is at most* $k(18L + 19\tau\Delta x)$

**Proof:** Note that when the procedure GoBackDown $(i, m, q)$ is invoked, $M_i = m, M_{i+1} = q, m \geq q + 2$, $X_i > X_{i+1}$, and the robot is at the post $P_i^m$. Let us consider the motion of the robot in each step of this procedure.

In step 1 the robot moves to the left, traversing edges associated with posts $P_j^r$ that are ahead of $P_{i+1}^q$ until a point $A$ is reached with $x$-coordinate $X_{i+1}^q$. (We associate an edge with the unique post on its right end.) Let $F_w$ be the highest fence reached in this step. Since no fence can have more than one post ahead of the last post of a lower fence (invariant [Almost-$x$-ordering]), the posts visited are precisely $P_i^m, P_{i-1}^m, \ldots, P_w^m$. For convenience, we say that these posts are "involved" in this procedure call.

In step 2 the robot goes greedy down-left until either (a) it is in a rectangular region whose opposite corners are the centers of the posts $(P_j^{m-1}, P_{j+1}^{m-1})$ where $w \leq j \leq i-1$, or (b) it is at $y$-coordinate $Y_i^{m-1}$. One of (a) or (b) must occur since the rectangular regions extend from the $y$-coordinate of $A$ down to $Y_i^{m-1}$. In the worst case, (a) occurs; in particular, the robot moves horizontally to the left and hits post $P_w^{m-1}$, after which it executes a sequence of calls GoDown $(j, m-1, m-1)$ for $j = w, w+1, \ldots, i-1$. We know from the proof of the previous Lemma 7 that after these calls the robot will be at post $P_i^{m-1}$. From that lemma, we can bound the cost of a call GoDown $(j, m-1, m-1)$ by $\tau$ plus twice the length of those edges associated with fence $F_{j+1}$ that lie between $x = X_j^{m-1}$ and $x = X_{j+1}^{m-1}$.

In step 3 the robot starts at some point $(x, y)$ with $y = Y_i^{m-1}$ and $X_i^{m-1} \leq x \leq X_{i+1}^q$ – this satisfies the conditions for GoDown $(i, m-1, q)$ and by the previous Lemma 7, the robot will arrive at $P_{i+1}^q$. This proves (1).

To bound the total cost of all calls to GoBackDown we show that after a call GoBackDown $(i, m, q)$ is made, (a) no future GoBackDown call can "involve" (in the above sense) the posts involved in this call, and (b) the *next* call GoBackDown $(i, m', q')$ from the *same* fence $F_i$ must have $q' \geq m$. Fact (a) implies that edges traversed by the robot in step 1 in two calls to GoBackDown are distinct; and the same holds for step 2. Thus the total step 1 cost is at most $k(3L + 3\tau\Delta x)$, which is the bound on the total length of the tree edges (Theorem 2). Also, the total step 2 cost (the cost of the GoDown calls) is at most $k\tau\Delta x$ plus *twice* the total edge-length, or $k(6L + 7\tau\Delta x)$. As in Lemma 7, fact (b) implies that the total step 3 cost is at most $k(9L + 9\tau\Delta x)$. Thus the contribution of GoBackDown to the total cost of FindFenceTree is at most $k(18L + 19\tau\Delta x)$, which proves (2).

We now argue facts (a) and (b). If a post $P_j^m$ is "involved" in a call to GoBackDown, it must be the case that $P_j^m$ is ahead of some lower fence that has no more than $m-2$ posts. Clearly, if $P_j^m$ is involved in the call GoBackDown $(i, m, q)$, then it cannot be involved again before the robot returns to $F_i$ since there is as yet no down-edge from $P_i^m$. On the other hand, if the robot returns to fence $F_i$, *every* fence below $F_i$ that is behind $F_i$ will have $m-1$ posts (Lemma 6). Therefore every fence that is below and behind $P_j^m$ will have $m-1$ posts, so post $P_j^m$ cannot be involved again. This proves fact (a). Fact (b) follows from Lemma 6: when the robot returns to $F_i$, after the call GoBackDown $(i, m, q)$, $M_{i+1}$ is either $m$ or $m-1$. If $M_{i+1} = m-1$ then case 2 of procedure FindFenceTree would apply and the robot would add the $m$th post to $F_{i+1}$. So in the next call GoBackDown $(i, m', q')$, $q'$ must be at least $m$. ■

We are thus able to bound the total cost of procedure FindFenceTree.

**Theorem 9** *The total distance walked by the robot in procedure* FindFenceTree *is at most* $k(60L + 62\tau\Delta x)$.

**Proof:** The motion of the robot is of 4 types: (a) finding a new edge, in cases 2 and 3, (b) going down to the next lower fence using procedure GoDown, in case 4, (c) going to the next lower fence using GoBackDown in case 1, and (d) returning to the fence above by following the path that was last taken to get to the current fence (case 6). Note that the total cost of (a) is the total length of all tree edges, and that the total cost of (d) is at most the sum of the costs of (a), (b), (c). The theorem then follows from Theorem 2, Lemma 7 and Lemma 8. ■